\documentclass[12pt]{article}
\usepackage{amssymb,latexsym,amsmath}
\usepackage{amsfonts}
\usepackage[cp1251]{inputenc}
\usepackage[epsf]{}
\usepackage[english,russian]{babel}
\usepackage[dvips]{graphicx}
\newtheorem{thm}{Theorem}
\newtheorem{prop}{Proposition}
\newtheorem{lm}{Lemma}
\newtheorem{cor}{Corollary}
\begin{document}
\title{Smale regular and chaotic A-homeomorphisms and A-diffeomorphisms}
\author{Medvedev V.$^{*}$\and Zhuzhoma E.\footnote{National Research University Higher School of Economics, 25/12 Bolshaya Pecherskaya, 603005, Nizhni Novgorod, Russia}\,\,\footnote{Corresponding author: zhuzhoma@mail.ru}} % }
\date{Dedicated to the memory of Stepin A.M.}
\maketitle

\renewcommand{\figurename}{Figure}
\renewcommand{\abstractname}{Abstract}
\renewcommand{\refname}{Bibliography}

\begin{abstract}
We introduce Smale A-homeomorphisms that includes regular, semi-chaotic, chaotic, and super chaotic homeomorphisms of topo\-lo\-gi\-cal $n$-manifold $M^n$, $n\geq 2$. Smale A-homeo\-mor\-p\-hisms contain A-diffeomorphisms (in particular, structurally stable diffeomorphisms) provided $M^n$ admits a smooth structure. Regular A-homeomorphisms contain all Morse-Smale diffeomorphisms, while semi-chaotic and chaotic A-homeomorphisms contain A-diffeo\-mor\-p\-hisms with trivial and nontrivial basic sets. Super chaotic A-homeo\-mor\-p\-hisms contain A-diffeomorphisms whose basic sets are nontrivial. We describe invariant sets that determine completely dynamics of regular, semi-chaotic, and chaotic Smale A-homeo\-mor\-p\-hisms. This allows us to get necessary and sufficient conditions of conjugacy for these Smale A-homeomorp\-hisms. We apply this necessary and sufficient conditions for structu\-ral\-ly stable surface diffeomorphisms with arbitrary number of one-dimensional expanding attractors. We also use this conditions to get the complete classification of Morse-Smale diffeomorphisms on projective-like $n$-manifolds for $n=2,8,16$.
\end{abstract}



\section*{Introduction}
Diffeomorphisms satisfying Smale's axiom A (in short, A-diffeomorphisms) were introduced by Smale \cite{Smale67} as a magnificent and natural generalization of structurally stable diffeomorphisms. By definition, a non-wandering set of A-diffeomorphism has a uniform hyperbolic struc\-tu\-re and is the topological closure of periodic orbits. Smale proved that the non-wandering set splits into closed, transitive, and invariant pieces called basic sets. A basic set is \textit{trivial}, if it is an isolated periodic orbit. A good example of A-diffeomorphism with trivial basic sets is a Morse-Smale diffeomorphism \cite{Palis69,Smale60a}. Such diffeomorphisms demonstrate regular dynamics. Due to Bowen \cite{Bowen1970b}, A-diffeomorphisms with nontri\-vial basic sets demonstrate chaotic dynamics since any such diffeomorphism has a positive entropy. The most familiar nontrivial basic sets are Plykin's attractor \cite{Plykin84} and codimension one expanding attractors introduced by Williams \cite{Williams1970a,Williams1974}. Such basic sets appeared in various applications, see for example \cite{GrinesMedvedevPochZh2015a,KuznetsovSP-2005,Strogatz-book-1994}.

Taking in mind that there are manifolds that do not admit smooth structures \cite{Milnor-1956}, we introduce Smale A-homeomorphisms with non-wandering sets having a hyperbolic type (see a precise definition below). Such homeomorphisms naturally appear in topological dynamical systems. For example, in \cite{EellsKuiper62}, it was proved the existence of topological Morse functions with three critical points on topological (including nonsmoothable) closed manifolds. Starting with these examples, one can construct topological (maybe, only topological) Morse-Smale flows and Morse-Smale homeomorphisms with the non-wandering set consisting of three fixed points of hyperbolic type. Deep theory of topological dynamical systems was developed in \cite{Akin-book-1993,AkinHurleyKennedy-book-2003}.

The challenging problem in the Theory of Dynamical Systems is the classification up to conjugacy dynamical systems with regular and chaotic dynamics.
Recall that homeo\-mor\-p\-hisms $f_1$, $f_2: M^n\to M^n$ are called \textit{conjugate}, if there is a homeomorp\-hism $h: M^n\to M^n$ such that $h\circ f_1=f_2\circ h$.
To check whether given $f_1$ and $f_2$ are conjugate, one constructs usually an invariant of conjugacy which is a dynamical characteristic keeping under a conjugacy homeomorphism. Normally, such invariant is constructed in the frame of special class of dynamical systems. The famous invariant is Poincare's rotation number for the class of transitive circle homeomorphisms \cite{Po1886}. This invariant is effective i.e. two transitive circle homeomorphisms are conjugate if and only if they have the same Poincare's rotation number (see \cite{NikZ99} and \cite{AnosovZhuzhoma2005}, ch. 7, concerning invariants of low dimensional dynamical systems).
Anosov \cite{Anosov67} and Smale \cite{Smale67} were first who realize the fundamental role of hyperbolicity for dynamical systems. Numerous topological invariants were constructed for special classes of A-diffeomorphisms including Anosov systems \cite{Fr70,Manning74,Newhouse70} and Morse-Smale systems, see the books \cite{ABZ,GrinMedvPo-book} and the surveys \cite{GrinesGurevichZhPochinka2019,MedvedevZhuzhoma2008-MIAN}.

In the frame of Smale A-homeomorphisms, we introduce regular, semi-chaotic, chaotic, and super chaotic homeomorphisms. We get necessary and sufficient conditions of conjugacy for regular, semi-chaotic, and chaotic Smale A-homeomorphisms on a closed topo\-lo\-gi\-cal $n$-manifold $M^n$, $n\geq 2$. Automatically, this gives necessary and sufficient conditions of conjugacy for Morse-Smale diffeomorphisms and a wide class of A-diffeomorphisms with nontrivial basic sets provided $M^n$ admits a smooth structure. We apply our conditions for structu\-ral\-ly stable surface diffeomor\-p\-hisms with arbitrary number of one-dimensional expanding attractors. We classify Morse-Smale dif\-feo\-morphisms with three periodic points on high-dimensional projective-like manifolds. Note that the projective-like manifolds were introduced by the authors in \cite{MedvedevZhuzhoma2016} (see also \cite{MedvedevZhuzhoma2013-top-appl}).

Let us give the main definitions and formulate the main results. Later on, $clos\,N$ means the topological closure of $N$. In \cite{MedvedevZhuzhoma2016}, the authors introduced the notation of equivalent embedding as follows.
Let $M^k_1$, $M^k_2\subset M^n$ be topologically embedded $k$-manifolds, $1\leq k\leq n-1$. We say they \textit{have the equivalent embedding} if there are neighborhoods $U(clos\, M^k_1)$, $U(clos\, M^k_2)$ of $clos\, M^k_1$, $clos\, M^k_2$ respectively and a homeomorphism $h: U(clos\, M^k_1)\to U(clos\, M^k_2)$ such that $h(M^k_1)=M^k_2$. This notation allows to classify Morse-Smale topological flows with non-wandering sets consisting of three equilibriums \cite{MedvedevZhuzhoma2016}. To be precise, it was proved that two such flows $f^t_1$, $f^t_2$ are topologically equivalent if and only if the stable (or unstable) separatrices of saddles of $f^t_1$, $f^t_2$ have the equivalent embedding. Remark that the notation of equivalent embedding goes back to a scheme introduced by Leontovich and Maier \cite{Leont-Maier37,Leont-Maier55} to attack the classification problem for flows on 2-sphere.

Solving the conjugacy problem for homeomorphisms, we have to add conjugacy relations to the equivalent embedding.
The modification of (global) conjugacy is a local conjugacy when the conjugacy holds in some neighborhoods of compact invariant sets. We introduce the intermediate notion, so-called a locally equivalent dynamical embedding (in short, dynamical embedding), as follows.
Let $f_1$, $f_2: M^n\to M^n$ be homeomorphisms of closed topological $n$-manifold $M^n$, $n\geq 2$, and $N_1$, $N_2$ invariant sets of $f_1$, $f_2$ respectively i.e. $f_i(N_i)=N_i$, $i=1,2$. We say that the sets $N_1$, $N_2$ \textit{have the same dynamical embedding} if there are neighborhoods $\delta_1$, $\delta_2$ of $clos~N_1$, $clos~N_2$ respectively and a homeomorphism $h_0: \delta_1\cup f_1(\delta_1)\to M^n$ such that
\begin{equation}\label{eq:dynamical-local-equivalent-embedding}
    h_0(\delta_1)=\delta_2,\qquad h_0(clos~N_1)=clos~N_2,\qquad h_0\circ f_1|_{\delta_1}=f_2\circ h_0|_{\delta_1}
\end{equation}

Recall that $F: L^n\to L^n$ is an A-diffeomorphism of smooth manifold $L^n$ provided the non-wandering set $NW(F)$ is hyperbolic, and the periodic orbits of $F$ are dense in $NW(F)$ \cite{Smale67}. The hyperbolicity implies that every point $z_0\in NW(F)$ has the stable $W^s(z_0)$ and unstable $W^u(z_0)$ manifolds formed by points $y\in L^n$ such that $\varrho_L(F^{pk}z_0, F^{pk}y)\to 0$ as
$k\to +\infty$ and $k\to -\infty$ respectively, where $\varrho _L$ is a metric on $L^n$ \cite{Grobman1959,Grobman1962,Hartman1963,HirschPughShub77-book,Pugh1962,Smale67}. Moreover, $W^s(z_0)$ and $W^u(z_0)$ are homeomorphic (in the interior topology) to Euclidean spaces $\mathbb{R}^{\dim W^s(z_0)}$, $\mathbb{R}^{\dim W^u(z_0)}$ respectively. Note that $\dim W^s(z_0) + \dim W^u(z_0)=n$. The non-wandering set $NW(F)$ is a finite union of pairwise disjoint $F$-invariant closed sets $\Omega _1$, $\ldots , \Omega _k$ such that every restriction $F|_{\Omega _i}$ is topologically transitive. These $\Omega _i$ are called \textit{basic sets} of $F$. A basic set is \textit{nontrivial} if it is not a periodic isolated orbit. Set $W^{s(u)}(\Omega_i)=\cup_{x\in\Omega_i}W^{s(u)}(x)$. One says that $\Omega _i$ is a \textit{sink (source) basic set} provided $W^u(\Omega _i)=\Omega _i$ ($W^s(\Omega _i)=\Omega _i$). A basic set $\Omega _i$ is a \textit{saddle basic set} if it is neither a sink nor a source basic set.

A homeomorphism $f: M^n\to M^n$ is called a \textit{Smale A-homeomorphism} if there is an A-diffeomorphism $F: L^n\to L^n$ such that the non-wandering sets $NW(f)$, $NW(F)$ have the same dynamical embedding. As a consequence, $NW(f)$ is a finite union of pairwise disjoint $f$-invariant closed sets $\Lambda _1$, $\ldots , \Lambda _k$ called \textit{basic sets} of $f$ such that every restriction $f|_{\Lambda _i}$ is topologically transitive. Each basic set $\Lambda$ has the stable manifold $W^s(\Lambda)$, and the unstable manifold $W^u(\Lambda)$. Similarly, one introduces the families of sink basic sets $\omega(f)$, and source basic sets $\alpha(f)$, and saddle basic sets $\sigma(f)$.

A Smale A-homeomorphism $f$ is called \textit{regular} if all basic sets $\omega(f)$, $\sigma(f)$, $\alpha(f)$ are trivial.

A Smale A-homeomorphism $f$ is called \textit{semi-chaotic} if exactly one family from the families $\omega(f)$, $\sigma(f)$, $\alpha(f)$ consists of non-trivial basic sets.

A Smale A-homeomorphism $f$ is called \textit{chaotic} if exactly two families from the families $\omega(f)$, $\sigma(f)$, $\alpha(f)$ consists of non-trivial basic sets.

A Smale A-homeomorphism $f$ is called \textit{super chaotic} if the families $\omega(f)$, $\sigma(f)$, $\alpha(f)$ consists of non-trivial basic  sets.

% Denote by $SsH(M^n)$ the set of Smale semi-chaotic homeomorphisms $M^n\to M^n$.

In Section \ref{s:exm}, we represent examples of all types above of Smale A-homeomorphisms. Actually, all examples are A-diffeomorphisms.

Now let us introduce invariant sets that determine dynamics of Smale homeo\-mor\-p\-hisms. Given any Smale A-homeomorphism $f: M^n\to M^n$, denote by $A(f)$ (resp., $R(f)$) the union of $\omega(f)$ (resp., $\alpha(f)$) and unstable (resp., stable) manifolds of saddle basic sets $\sigma(f)$ :
 $$ A(f)=\omega(f)\bigcup_{\nu\in\sigma(f)}W^u(\nu),\quad R(f)=\alpha(f)\bigcup_{\nu\in\sigma(f)}W^s(\nu). $$

The following statement gives the necessary and sufficient conditions of conjugacy for three types of Smale A-homeomorphisms. Note that if $f$ is a regular or semi-chaotic Smale A-homeomor\-p\-hism, then at least one of the families $\omega(f)$, $\alpha(f)$ consists of trivial basic sets. However, if $f$ is a chaotic Smale A-homeomorphism, then it is possible a priori that the both $\omega(f)$ and $\alpha(f)$ consist of nontrivial basic sets but $\sigma(f)$ consists of trivial basic sets.
\begin{thm}\label{thm:congugacy-semi-regular}
Let $M^n$ be a closed topological $n$-manifold $M^n$, $n\geq 2$ and $f_i: M^n\to M^n$ is either regular, or semi-chaotic, or chaotic Smale A-homeomorphism, $i=1,2$. For chaotic $f_i$, we suppose that either $\omega(f_i)$ or $\alpha(f_i)$ consists of trivial basic sets, $i=1,2$. Then the homeomorphisms $f_1$, $f_2$ are conjugate if and only if one of the following conditions holds:
\begin{itemize}
  \item the basic sets $\alpha(f_1)$, $\alpha(f_2)$ are trivial while the sets $A(f_1)$, $A(f_2)$ have the same dynamical embedding;
  \item the basic sets $\omega(f_1)$, $\omega(f_2)$ are trivial while the sets $R(f_1)$, $R(f_2)$ have the same dynamical embedding.
\end{itemize}
\end{thm}

In Section \ref{s:discus-and-appl}, we apply Theorem \ref{thm:congugacy-semi-regular} to consider the conjugacy for structurally stable surface diffeomorphisms $M^2\to M^2$ with one-dimensional (orientable and non-orientable) attractors $\Lambda_1$, $\ldots$, $\Lambda_k$, $k\geq 2$ (remark that a structurally stable diffeomorphism is an A-diffeomorphism \cite{Mane88b}). We also use Theorem \ref{thm:congugacy-semi-regular} to classify Morse-Smale dif\-feo\-morphisms with three periodic points on projective-like $n$-manifolds for $n\in\{2,8,16\}$.

\begin{figure}[h]
\centerline{\includegraphics[height=35mm]{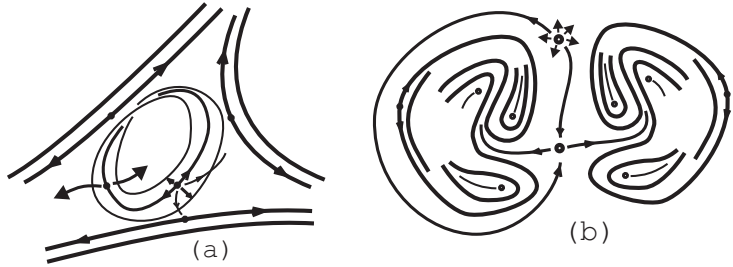}}
\caption{(a) one isolated saddle and one expanding attractor on some non-oriented surface; (b) one isolated saddle and two Plykin attractors.} \label{two-exp+non-or}
\end{figure}

First, we prove the following statement interesting itself (note that a one dimensional expanding attractor is a trivial basic set).
\begin{prop}\label{prop:attractors-no-more-saddles-plus-one}
Let $f: M^2\to M^2$ be an A-diffeomorphism with the non-wandering set $NW(f)$ consisting of one-dimensional expanding attractors $\Lambda_1$, $\ldots$, $\Lambda_k$, and $s_0\geq 0$ isolated saddle periodic points, and arbitrary number of isolated nodal periodic orbits. Then $k\leq s_0+1$.
\end{prop}
The case $k=s_0=1$ is represented in Fig.~\ref{two-exp+non-or}, (a), while the case $k=2$ and $s_0=1$ is represented in Fig.~\ref{two-exp+non-or}, (b) with two Plykin attractors. See also \cite{AransonPlykinZhirovZhuzhoma1997} where one got the estimate for the number of one-dimensional basic sets of surface A-diffeomorphisms depending on a genus of supporting surface.

The following statement shows that the dynamical embedding of unstable manifolds of isolated saddles (trivial basic sets) determine completely global dynamics of structurally stable surface diffeomorphisms with one-dimensional expanding attractors.
\begin{thm}\label{thm:many-expand-attr-one-less-saddle}
Let $f_i: M^2\to M^2$, $i=1,2$, be a structurally stable diffeomorphism of closed 2-manifold $M^2$ such that the spectral decomposition of $f_i$ consists of $k\geq 2$ one-dimensional expanding attractors $\Lambda_1^{(i)}$, $\ldots$, $\Lambda_k^{(i)}$, and isolated source periodic orbits, and $k-1$ isolated saddle periodic points denoted by $\sigma_1^{(i)}$, $\ldots$, $\sigma_{k-1}^{(i)}$. Then $f_1$, $f_2$ are conjugate if and only if the sets $\cup_{j=1}^{j=k-1}W^u(\sigma_j^{(1)})$, $\cup_{j=1}^{j=k-1}W^u(\sigma_j^{(2)})$ have the same dynamical embedding.
\end{thm}

Denote by $SRH(M^n)$ the class of Smale regular homeomorphisms $M^n\to M^n$. Note that it is possible that $f\in SRH(M^n)$ has the empty set $\sigma(f)$ of saddle periodic points. In this case the set $\alpha(f)$ consists of a unique source and the set $\omega(f)$ consists of a unique sink, and $M^n=S^n$ is an $n$-sphere. Later on, we'll assume that $f\in SRH(M^n)$ has a non-empty set $\sigma(f)$ of saddle periodic points.

%\begin{figure}[h]
%\centerline{\includegraphics[height=3.5cm]{possible.eps}}\caption{Examples of regular Smale diffeomorphisms.} \label{possible}
%\end{figure}

Clearly, $SRH(M^n)$ contains all Morse-Smale diffeomorphisms provided $M^n$ admits a smooth structure. Note that the class $SRH(M^n)$ is an essential extension of Morse-Smale diffeomorphisms for a Smale regular
diffeomorphism can contain nonhyperbolic periodic points, tangencies, and separatrix connections.
%Fig.~\ref{possible}.

As a consequence of Theorem \ref{thm:congugacy-semi-regular}, one gets the following statement (in particular, one gets the necessary and sufficient conditions of conjugacy for any Morse-Smale diffeomorphisms on smooth closed manifolds).
\begin{cor}\label{cor:congugacy-regular}
Let $M^n$ be a closed topological $n$-manifold $M^n$, $n\geq 2$. Homeomorphisms $f_1$, $f_2\in SRH(M^n)$ are conjugate if and only if one of the following conditions holds:
\begin{itemize}
  \item the sets $A(f_1)$, $A(f_2)$ have the same dynamical embedding;
  \item the sets $R(f_1)$, $R(f_2)$ have the same dynamical embedding.
\end{itemize}
\end{cor}

Denote by $MS(M^n;a,b,c)$ the set of Morse-Smale diffeomorphisms $f: M^n\to M^n$ whose non-wandering set consists of $a$ sinks, $b$ sources, and $c$ saddles. In \cite{MedvedevZhuzhoma2013-top-appl}, the authors proved that for $MS(M^n;1,1,1)$ the only values of $n$ possible are $n\in\{2,4,8,16\}$. Moreover, the supporting manifolds for $MS(M^n;1,1,1)$ are projective-like provided $n\in\{2,8,16\}$  \cite{MedvedevZhuzhoma2013-top-appl,MedvedevZhuzhoma2016}.

First, to illustrate the applicability of Corollary \ref{cor:congugacy-regular}, we consider very simple class $MS(M^2;1,1,1)$. In this case, the supporting manifold $M^2$ is the projective plane $M^2=\mathbb{P}^2$ \cite{MedvedevZhuzhoma2013-top-appl}. Below, we define a type for a unique saddle of $f\in MS(\mathbb{P}^2,1,1,1)$. Using Corollary \ref{cor:congugacy-regular} we'll show how to get the following complete classification of Morse-Smale diffeomorphisms $MS(\mathbb{P}^2,1,1,1)$.
\begin{prop}\label{prop:proj-plane-top-conj}
Two diffeomorphisms $f_1$, $f_2\in MS(\mathbb{P}^2,1,1,1)$ are conjugate if and only if the types of their saddles coincide. There are four types $T_i$, $i=1,2,3,4$ of a saddle. Given any type $T_i$, $i=1,2,3,4$, there is a diffeomorphism $f\in MS(\mathbb{P}^2,1,1,1)$ with a saddle $\sigma(f)$ of the type $T_i$.
\end{prop}
Thus, up to conjugacy, there are four classes of Morse-Smale diffeomorphisms $MS(\mathbb{P}^2,1,1,1)$.

The most essential application is a complete classification of Morse-Smale diffeomorphisms $MS(M^8;1,1,1)$ and $MS(M^{16};1,1,1)$. The supporting $2k$-manifolds for diffeomorphisms from the set $MS^{2k}(1,1,1)$ will be denoted by $M^{2k}(1,1,1)$. Recall that $S^k$ is a $k$-sphere. Below, $\alpha_f$, $\sigma_f$, and $\omega_f$ mean the source, the saddle, and the sink of
$f\in MS^{2k}(1,1,1)$ respectively.

An embedding $\varphi: S^k\to M^{2k}(1,1,1)$ is called \textit{basic} if
\begin{itemize}
  \item $\varphi(S^k)$ is a locally flat $k$-sphere;
  \item $M^{2k}(1,1,1)\setminus \varphi(S^k)$ is an open $2k$-ball, $M^{2k}(1,1,1)=B^{2k}\sqcup\varphi(S^k)$.
\end{itemize}
It was proved in \cite{MedvedevZhuzhoma2013-top-appl} that every supporting manifold $M^{2k}(1,1,1)$, $k=4,8$, admits a basic embedding.

\begin{thm}\label{thm:classif-M-2k-for-8-16}
Let $f: M^{2k}(1,1,1)\to M^{2k}(1,1,1)$ be a diffeomorphism from the set $MS^{2k}(1,1,1)$, $k=4,8$. Then the following claims hold :

1) for any $f\in MS^{2k}(1,1,1)$, there are basic embedding
 $$\varphi_u(f): S^k\to M^{2k}(1,1,1),\qquad \varphi_s(f): S^k\to M^{2k}(1,1,1)$$
such that $\varphi_u(f)(S^k)=W^u_{\sigma_f}\cup\{\omega_f\}$ and $\varphi_s(f)(S^k)=W^s_{\sigma_f}\cup\{\alpha_f\}$.

2) given any basic embedding $\varphi: S^k\to M^{2k}(1,1,1)$, there is $f\in MS^{2k}(1,1,1)$ such that one of the following equality holds:
 $$\varphi(S^k)=W^u_{\sigma_f}\cup\{\omega_f\}\,\, \mbox{ or }\,\, \varphi(S^k)=W^s_{\sigma_f}\cup\{\alpha_f\} $$

3) two Morse-Smale diffeomorphisms $f_1, f_2\in MS^{2k}(1,1,1)$ are conjugate if and only if one of the following conditions holds:
\begin{itemize}
  \item the basic embedding $\varphi_u(f_1)(S^k)=W^u_{\sigma_{f_1}}\cup\{\omega_{f_1}\}$, $\varphi_u(f_2)(S^k)=W^u_{\sigma_{f_2}}\cup\{\omega_{f_2}\}$ have the same dynamical embedding;
  \item the basic embedding $\varphi_s(f_1)(S^k)=W^s_{\sigma_{f_1}}\cup\{\alpha_{f_1}\}$, $\varphi_s(f_2)(S^k)=W^s_{\sigma_{f_2}}\cup\{\alpha_{f_2}\}$ have the same dynamical embedding.
\end{itemize}
\end{thm}
Thus, every $f\in MS^{2k}(1,1,1)$ corresponds the basic embedding $\varphi(f): S^k\to M^{2k}(1,1,1)$. Given any basic embedding $\varphi$, there is $f\in MS^{2k}(1,1,1)$ such that $\varphi(f)=\varphi$. At last, a dynamical embedding of basic embedding defines completely a conjugacy class in $MS^{2k}(1,1,1)$. We see, that the set of basic embedding (up to isotopy) form the admissible set of conjugacy invariants for the Morse-Smale diffeomorphisms $MS^{2k}(1,1,1)$, $k=4,8$.
As to the class $MS(M^4;1,1,1)$, the existence of realizable and effective conjugacy invariant is still the open problem. The main reason is the possibility of wild embedding for topological closure of separatrix of a saddle \cite{MedvedevZhuzhoma2013-top-appl}.

The next statement shows that an equivalent embedding gives an invariant of conjugacy for iterations of Morse-Smale diffeomorphisms $f\in MS^{2k}(1,1,1)$, $k=4,8$.
\begin{thm}\label{thm:iterations-equiv-embed}
Let $\sigma_{f_i}$ be a unique saddle of $f_i\in MS^{2k}(1,1,1)$, $i=1,2$. If the stable (unstable) manifolds $W^{s(u)}(\sigma_{f_1})$, $W^{s(u)}(\sigma_{f_2})$ have equivalent embedding,
then there are $k_1$, $k_2\in\mathbb{N}$ such that the diffeomorphisms $f_1^{k_1}$, $f_2^{k_2}$ are conjugate.
Vice versa, if $f_1^{k_1}$, $f_2^{k_2}$ are conjugate for some $k_1$, $k_2\in\mathbb{N}$, then the stable (unstable) manifolds $W^{s(u)}(\sigma_{f_1})$, $W^{s(u)}(\sigma_{f_2})$ have equivalent embedding.
\end{thm}

The structure of the paper is the following. In Section \ref{s:def-and-prev}, we give some previous results. In Section \ref{s:proof-hom}, we prove Theorem \ref{thm:congugacy-semi-regular}.
In Section \ref{s:discus-and-appl}, we prove Propositions \ref{prop:attractors-no-more-saddles-plus-one}, \ref{prop:proj-plane-top-conj}, and Theorems \ref{thm:many-expand-attr-one-less-saddle}, \ref{thm:classif-M-2k-for-8-16}, \ref{thm:iterations-equiv-embed}.

\textsl{Acknowledgments}. This work is supported by Laboratory of Dynamical Systems and Appli\-ca\-tions of National Research University Higher School of Economics, of the Ministry of science and higher education of the RF, grant ag. № 075-15-2019-1931.

\section{Examples of A-diffeomorphisms}\label{s:exm}

1) \textsl{Regular A-diffeomorphisms}. Obvious example of regular A-diffeomorphism is a Morse-Smale diffeomorphism. Note that there are regular A-diffeomorphisms that do not belong to the set of Morse-Smale diffeomorphisms. For example, they can belong to the boundary of the set of Morse-Smale diffeomorphisms in the space of diffeomorphisms. There are regular A-diffeomorphisms which can not be approximated by Morse-Smale diffeomorphisms \cite{PughWalkerWilson1977}.

\begin{figure}[h]
\centerline{\includegraphics[height=35mm]{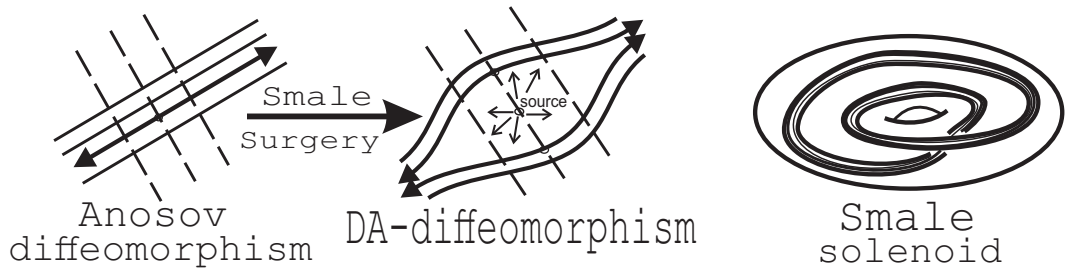}}\caption{Examples of semi-chaotic Smale diffeomorphisms.} \label{DA+SmaleSol}
\end{figure}

2) \textsl{Semi-chaotic A-diffeomorphisms}. A good example of semi-chaotic diffeomorphism is a so-called DA-diffeomorphism obtained from Anosov automorphism after Smale surgery \cite{Smale67}, see Fig.~\ref{DA+SmaleSol}.
A classical DA-diffeomorphism $f: \mathbb{T}^2\to\mathbb{T}^2$ contains a nontrivial attractor $\omega(f)$, trivial repeller $\alpha(f)$, and empty set $\sigma(f)$. A generalized DA-diffeomorphism contains nonempty set $\sigma(f)$ \cite{GinesZhuzhoma2005}. Another example of semi-chaotic A-diffeomorphism is a classical Smale horseshoe $g_s: \mathbb{S}^2\to\mathbb{S}^2$. Well-known that there is $g_s$ with trivial attractor $\omega(g_s)$ and repeller $\alpha(g_s)$, and nontrivial $\sigma(g_s)$.

Starting with DA-diffeomorphisms, Williams \cite{Williams1970a} constructed an open domain $\mathcal{N}\subset Diff^1(\mathbb{T}^2)$ consisting of structurally unstable diffeo\-mor\-p\-hisms. It is easy to see that $\mathcal{N}$ contains semi-chaotic A-diffeomorphisms.

\begin{figure}[h]
\centerline{\includegraphics[height=3.5cm]{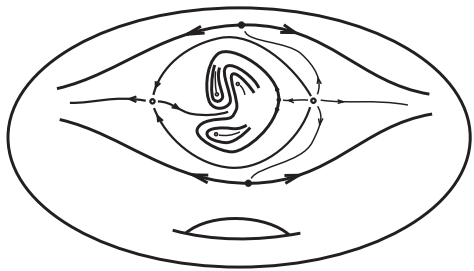}}\caption{One isolated saddle and two expanding attractors on a torus.} \label{plyk+DA}
\end{figure}

One more example of semi-chaotic A-diffeomorphism is shown in Fig.~\ref{plyk+DA} with a DA-attractor and Plykin attractor on a torus.

3) \textsl{Chaotic A-diffeomorphisms}. Take the classical DA-diffeomorphism $f: \mathbb{T}^2\to\mathbb{T}^2$ with the non-wandering set consisting of a source $\alpha$ and one-dimensional expanding attractor $\Lambda_a$. The diffeomorphism $f^{-1}$ defined on a copy $\mathbb{T}^2$ has the non-wandering set consisting of a sink $\omega$ and one-dimensional contracting repeller $\Lambda_r$. Let us delete a small neighborhood $U_a$ (resp., $U_s$) of $\alpha$ (resp., $\omega$) homeomorphic to a disk. Take an orientation reversing diffeomorphism $h: \partial U_a\to\partial U_r$. Then the surface
$M^2=\left(\mathbb{T}^2\setminus U_a\right)\cup_h\left(\mathbb{T}^2\setminus U_r\right)$ is a pretzel (closed orientable surface of genus 2). Following \cite{RobinsonWilliams71}, one can construct A-diffeomorphism
$g: M^2 \to M^2$ with the non-wandering set consisting of $\Lambda_a\cup\Lambda_r$ such that $g|_{\Lambda_a}=f$ and $g|_{\Lambda_r}=f^{-1}$. Thus, $\alpha(g)=\Lambda_r$ and $\omega(g)=\Lambda_a$. Clearly, $g$ is a chaotic A-diffeomorphism. Due to \cite{RobinsonWilliams71}, there is a construction such that $g$ has a closed simple curve consisting of the tangencies of the invariant stable manifolds of $\Lambda_a$ and the invariant unstable manifolds of $\Lambda_r$.

Another example one gets starting with a Smale solenoid \cite{Smale67}, see Fig.~\ref{DA+SmaleSol}. This mapping can be extended to $\Omega$-stable diffeomorphism $f_s: M^3\to M^3$ with a one-dimensional expanding attractor, say $\Omega_1$, and one-dimensional contracting repeller, say $\Omega_2$, where $M^3$ is a 3-sphere or lens space \cite{JiagNiWang2004,JimingMaBinYu2007}. This chaotic diffeomorphism is similar to the Robinson-Williams diffeomorphism $g$ considered above.
There is a bifurcation of $\Omega_1$ into a zero-dimensional saddle type basic set and isolated attracting periodic orbits \cite{ZhIsaenkova2011-transl}. As a result, one gets a chaotic Smale diffeomorphism $f_0: M^3\to M^3$ with trivial basic sets $\omega(f_0)$, and the nontrivial source basic set $\alpha(f_0)=\Omega_2$, and the nontrivial zero-dimensional saddle basic set $\sigma(f_0)$.

4) \textsl{Super chaotic A-diffeomorphisms}. Let $g_s: \mathbb{S}^2\to\mathbb{S}^2$ be the classical Smale horseshoe and $f: \mathbb{T}^2\to\mathbb{T}^2$ the classical DA-diffeomorphism considered above. Delete small neighborhoods $U_1$, $U_2$ of the sink $\omega(g_s)$ and the source $\alpha(g_s)$ respectively each homeomorphic to a disk. There are reversing orientation diffeomorphisms $h_1: \partial U_1\to\partial U_a$ and
$h_2: \partial U_2\to\partial U_r$. Then the surface
$M^2=\left(\mathbb{T}^2\setminus U_a\right)\bigcup_{h_1}\left(S^2\setminus U_1\cup U_2\right)\bigcup_{h_2}\left(\mathbb{T}^2\setminus U_r\right)$ is a pretzel. Similarly to Robinson-Williams's method developed in \cite{RobinsonWilliams71}, one can construct a diffeomorphism $g_0: M^2\to M^2$ with $\alpha(g_0)=\Lambda_r$, $\omega(g_0)=\Lambda_a$, and $\sigma(g_0)$ homeomorphic to the Smale horseshoe $\sigma(g_s)$. Thus, $g_0$ is a super chaotic A-diffeomorphism. By similar way, one can get another examples starting with semi-chaotic A-diffeomorphisms.

\section{Properties of Smale homeomorphisms}\label{s:def-and-prev}\nopagebreak

We begin by recalling several definitions. Further details may be found in \cite{AnosovZhuzhoma2005,ABZ,Smale67}.
Denote by $Orb(x)$ the orbit of point $x\in M^n$ under a homeomorphism $f: M^n\to M^n$. The $\omega$-limit set $\omega(x)$ of the point $x$ consists of the points $y\in M^n$ such that $f^{k_i}(x)\to y$ for some sequence $k_i\to\infty$. Clearly that any points of $Orb(x)$ have the same $\omega$-limit. Replacing $f$ with $f^{-1}$, one gets an $\alpha$-limit set. Obviously, $\omega(x)\cup\alpha(x)\subset NW(f)$ for every $x\in M^n$.

Later on, $f\in SsH(M^n)$. Given a family $C=\{c_1,\ldots,c_l\}$ of sets $c_i\subset M^n$, denote by $\widetilde{C}$ the union $c_1\cup\ldots\cup c_l$. It follows immediately from definitions that
\begin{equation}\label{eq:non-wandering-for-AP}
    NW(f)=\widetilde{\alpha(f)}\cup\widetilde{\omega(f)}\cup\widetilde{\sigma(f)},\quad f\in SsH(M^n)
\end{equation}

\begin{lm}\label{lm:x-go-to-saddle} Let $f\in SsH(M^n)$ and $x\in M^n$. Then
\begin{enumerate}
  \item if $\omega(x)\subset\widetilde{\sigma(f)}$, then $x\in W^s(\sigma_*)$ for some saddle basic set $\sigma_*\in\sigma(f)$.
  \item if $\alpha(x)\subset\widetilde{\sigma(f)}$, then $x\in W^u(\sigma_*)$ for some saddle basic set $\sigma_*\in\sigma(f)$.
\end{enumerate}
\end{lm}
\textsl{Proof}. Suppose that $\omega(x)\subset\widetilde{\sigma(f)}$. Since $\widetilde{\alpha(f)}$ and $\widetilde{\omega(f)}$ are invariant sets, $x\notin\widetilde{\alpha(f)}\cup\widetilde{\omega(f)}$. Therefore, there are exist a neighborhood $U(\alpha)$ of $\alpha(f)$ and neighborhood $U(\omega)$ of $\omega(f)$ such that the positive semi-orbit $Orb^+(x)$ belongs to the compact set $N=M^n\setminus\left(U(\omega)\cup U(\alpha)\right)$. Let $V(\sigma_1)$, $\ldots$, $V(\sigma_m)$ be pairwise disjoint neighborhoods of saddle basic sets $\sigma_1$, $\ldots$, $\sigma_m$ respectively such that $\cup_{i=1}^mV(\sigma_i)\subset N$. Since every $V(\sigma_i)$ does not intersect $\cup_{j\neq i}V(\sigma_j)$ and all saddle basic sets are invariant, one can take the neighborhoods $V(\sigma_1)$, $\ldots$, $V(\sigma_m)$ so small that every $f(V(\sigma_i))$ does not intersect
$\cup_{j\neq i}V(\sigma_j)$. Suppose the contrary, i.e. there is no a unique saddle basic set $\sigma_*\in\sigma(f)$ with $x\in W^s(\sigma_*)$. Thus, there are at least two different saddle basic sets $\sigma_1$, $\sigma_2$ such that $x\in W^s(\sigma_1)$ and $x\in W^s(\sigma_2)$. Hence, $\omega(x)$ have to intersect $\sigma_1$, $\sigma_2$. It follows that the compact set $N_0=N\setminus\left(\cup_{i=1}^mV(\sigma_i)\right)$ contains infinitely many points of the semi-orbit $Orb^+(x)$. This implies $\omega(x)\cap N_0\neq\emptyset$ that contradicts (\ref{eq:non-wandering-for-AP}).
The second assertion is proved similarly.
$\Box$

A set $U$ is a \textit{trapping region} for $f$ if $ f\left(clos~U\right)\subset int~U.$ A set $A$ is an \textit{attracting set} for $f$ if there exists a trapping set $U$ such that
 $$ A=\bigcap_{k\geq 0}f^k(U). $$
A set $A^*$ is a \textit{repelling set} for $f$ if there exists a trapping region $U$ for $f$ such that
 $$ A^*=\bigcap_{k\leq 0}f^k(M^n\setminus U). $$
Another words, $A^*$ is an attracting set for $f^{-1}$ with the trapping region $M^n\setminus U$ for $f^{-1}$. When we wish to emphasize the dependence of an attracting set $A$ or a repelling set $A^*$ on the trapping region $U$ from which it arises, we denote it by $A_U$ or $A^*_U$ respectively.

Let $A$ be an attracting set for $f$. The \textit{basin} $B(A)$ of $A$ is the union of all open trapping regions $U$ for $f$ such that $A_U=A$. One can similarly define the notion of basin for a repelling set.

Let $N$ be an attracting or repelling set and $B(N)$ the basin of $N$. A closed set $G(N)\subset B(N)\setminus N$ is called a \textit{generating set} for the domain $B(N)\setminus N$ if
 $$ B(N)\setminus N = \cup_{k\in\mathbb{Z}}f^k\left(G(N)\right). $$
Moreover,

1) every orbit from $B(N)\setminus N$ intersects $G(N)$; 2) if an orbit from $B(N)\setminus N$ intersects the interior of $G(N)$, then this orbit intersects $G(N)$ at a unique point; 3) if an orbit from $B(N)\setminus N$ intersects the boundary of $G(N)$, then the intersection of this orbit with $G(N)$ consists of two points; 4) the boundary of $G(N)$ is the union of finitely many compact codimension one topological submanifolds.

\begin{lm}\label{lm:A-is-attractive-set-alpha-trivial} Let $f\in SsH(M^n)$.

1) Suppose that all basic sets $\alpha(f)$ are trivial. Then $\widetilde{\alpha(f)}$ is a repelling set while $A(f)$ is an attracting set with
 $$ B\left(\widetilde{\alpha(f)}\right)\setminus\widetilde{\alpha(f)} = B\left(A(f)\right)\setminus A(f). $$
Moreover,
\begin{itemize}
  \item there is a trapping region $T(\alpha)$ for $f^{-1}$ of the set $\widetilde{\alpha(f)}$ consisting of pairwise disjoint open $n$-balls
        $b_1$, $\ldots$, $b_r$ such that each $b_i$ contains a unique periodic point from $\alpha(f)$;
  \item the regions $B(\widetilde{\alpha(f)})\setminus\widetilde{\alpha(f)}$, $B(A(f))\setminus A(f)$ have the same generating set $G(\alpha)$ consisting of pairwise disjoint
        closed $n$-annuluses $a_1$, $\ldots$, $a_r$ such that $a_i=clos~f^{p_i}(b_i)\setminus b_i$ where $p_i\in\mathbb{N}$ is a minimal period of a periodic point belonging
        to $b_i$, $i=1,\ldots,r$ :
        $$ G(\alpha)=\cup_{i=1}^ra_i=\cup_{i=1}^r\left(clos~f^{p_i}(b_i)\setminus b_i\right); $$
  \item $B(A(f))\setminus A(f)=\cup_{k\in\mathbb{Z}}f^k(G(\alpha))$.
\end{itemize}

2) Suppose that all basic sets $\omega(f)$ are trivial. Then $\widetilde{\omega(f)}$ is an attracting set while $R(f)$ is a repelling set with
 $$ B(\widetilde{\omega(f)})\setminus\widetilde{\omega(f)} = B(R(f))\setminus R(f). $$
Moreover,
\begin{itemize}
  \item there is a trapping region $T(\omega)$ for $f$ of the set $\widetilde{\omega(f)}$ consisting of pairwise disjoint an open $n$-balls
        $b_1$, $\ldots$, $b_l$ such that each $b_i$ contains a unique periodic point from $\omega(f)$;
  \item the regions $B(\widetilde{\omega(f)})\setminus\widetilde{\omega(f)}$, $B(R(f))\setminus R(f)$ have the same generating set $G(\omega)$ consisting of pairwise disjoint
        closed $n$-annuluses $a_1$, $\ldots$, $a_l$ such that $a_i=b_i\setminus int~f^{p_i}(b_i)$ where $p_i\in\mathbb{N}$ is a minimal period of a periodic point belonging
        to $b_i$, $i=1,\ldots,l$ :
        $$ G(\omega)=\cup_{i=1}^ra_i=\cup_{i=1}^r\left(b_i\setminus int~f^{p_i}(b_i)\right); $$
  \item $B(R(f))\setminus R(f)=\cup_{k\in\mathbb{Z}}f^k(G(\omega))$.
\end{itemize}
\end{lm}
\textsl{Proof}. It is enough to prove the first statement only. Since all basic sets $\alpha(f)$ are trivial and consists of locally hyperbolic source periodic points, there is a trapping region $T(\alpha)$ for $f^{-1}$ of the set $\widetilde{\alpha(f)}$ consisting of pairwise disjoint open $n$-balls $b_1$, $\ldots$, $b_r$ such that each $b_i$ contains a unique periodic point $q_i$ from $\alpha(f)$ \cite{Palis69,Smale60a}. Thus,
 $$ T(\alpha)=\cup_{i=1}^rb_i,\quad \cap_{k\leq 0}f^{kp_i}(b_i)=q_i,\quad i=1,\ldots,r. $$
As a consequence, there is the generating set
$G(\alpha)=\cup_{i=1}^r\left(clos~f^{p_i}(b_i)\setminus b_i\right)$ consisting of pairwise disjoint closed $n$-annuluses $a_i=clos~f^{p_i}(b_i)\setminus b_i$, $i=1,\ldots,r$.

Since the balls $b_1$, $\ldots$, $b_r$ are pairwise disjoint and $clos~b_i\subset f^{p_i}(b_i)$, the balls $f^{p_1}(b_1)$, $\ldots$, $f^{p_r}(b_r)$ are pairwise disjoint also. For simplicity of exposition, we'll assume that $\alpha(f)$ consists of fixed points (otherwise, $\alpha(f)$ is divided into periodic orbits each considered like a point).
Therefore,
 $$ f\left(M^n\setminus\cup_{i=1}^rb_i\right)=M^n\setminus\cup_{i=1}^rf(b_i)\subset M^n\setminus\cup_{i=1}^rclos~b_i\subset int~\left(M^n\setminus\cup_{i=1}^rb_i\right). $$
Hence, $M^n\setminus\cup_{i=1}^rb_i$ is a trapping region for $f$. Clearly, $A(f)\subset M^n\setminus\cup_{i=1}^rb_i$.

Take a point $x\in M^n\setminus\cup_{i=1}^rb_i$. Obviously, $\omega(x)\notin\widetilde{\alpha(f)}$. It follows from (\ref{eq:non-wandering-for-AP}) that $\omega(x)\in\widetilde{\omega(f)}\cup\widetilde{\sigma(f)}$. By
Lemma \ref{lm:x-go-to-saddle}, $\omega(x)\in A(f)$. Therefore, $A(f)$ is an attracting set with the trapping region $M^n\setminus\cup_{i=1}^rb_i$ for $f$ :
 $$ A(f)=A_{M^n\setminus\cup_{i=1}^rb_i}. $$
Moreover,
 $$ M^n=\widetilde{\alpha(f)}\cup B(A(f)) $$
because of $\cap_{k\leq 0}f^k(b_i)=q_i$, $i=1,\ldots,r$.

Let us prove the quality $B\left(\widetilde{\alpha(f)}\right)\setminus\widetilde{\alpha(f)} = B\left(A(f)\right)\setminus A(f)$. Take
$x\in B\left(\widetilde{\alpha(f)}\right)\setminus\widetilde{\alpha(f)}$. Since $x\notin\widetilde{\alpha(f)}$ and $M^n=\widetilde{\alpha(f)}\cup B(A(f))$, $x\in B(A(f))$. Since
$x\in B\left(\widetilde{\alpha(f)}\right)$, $\alpha(x)\subset\alpha(f)$. Hence, $x\notin A(f)$ and $x\in B(A(f))\setminus A(f)$. Now, set $x\in B(A(f))\setminus A(f)$. Then $x\notin\alpha(f)$. Since $x\notin A(f)$, $\alpha(x)\subset \widetilde{\sigma(f)}\cup\widetilde{\alpha(f)}$. If one assumes that $\alpha(x)\subset\widetilde{\sigma(f)}$, then according to Lemma \ref{lm:x-go-to-saddle}, $x\in W^u(\nu)$ for some saddle basic set $\nu$. Thus, $x\in A(f)$ which contradicts to $x\notin A(f)$. Therefore, $\alpha(x)\subset \widetilde{\alpha(f)}$. Hence $x\in B\left(\widetilde{\alpha(f)}\right)$. As a consequence, $x\in B\left(\widetilde{\alpha(f)}\right)\setminus\widetilde{\alpha(f)}$.

The last assertion of the first statement follows from the previous ones. This completes the proof.
$\Box$

In the next statement, we keep the notation of Lemma \ref{lm:A-is-attractive-set-alpha-trivial}.
\begin{lm}\label{lm:neighbor-for-A-is-attractive-set-alpha-trivial} Let $f\in SsH(M^n)$.

1) Suppose that all basic sets $\alpha(f)$ are trivial. Then given any neighborhood $V_0(A)$ of $A(f)$, there is $n_0\in\mathbb{N}$ such that
 $$ \cup_{k\geq n_0}f^k\left(G(\alpha)\right)\subset V_0(A) $$
where $G(\alpha)$ is the generating set of the region $B(\widetilde{\alpha(f)})\setminus\widetilde{\alpha(f)}$.

2) Suppose that all basic sets $\omega(f)$ are trivial. Then given any neighborhood $V_0(R)$ of $R(f)$, there is $n_0\in\mathbb{N}$ such that
 $$ \cup_{k\leq -n_0}f^k\left(G(\omega)\right)\subset V_0(R) $$
where $G(\omega)$ is the generating set of the region $B(\widetilde{\omega(f)})\setminus\widetilde{\omega(f)}$.
\end{lm}
\textsl{Proof}. It is enough to prove the first statement only. Take a closed tripping neighborhood $U$ of $A(f)$ for $f$. Since $\cap_{k\in\mathbb{N}}f^k(U)=A(f)\subset V_0(A)$, there is $k_0\in\mathbb{N}$ such that $f^{k_0}(U)\subset V_0(A)$. Clearly, $f^{k_0}(U)$ is a tripping region of $A(f)$ for $f$. Hence, $f^{k_0+k}(U)\subset f^{k_0}(U)\subset V_0(A)$ for every $k\in\mathbb{N}$.

Let $G(\alpha)$ be a generating set of the region $B(\widetilde{\alpha(f)})\setminus\widetilde{\alpha(f)}$. By Lemma \ref{lm:A-is-attractive-set-alpha-trivial}, $G(\alpha)$ is the generating set of the region $B\left(A(f)\right)\setminus A(f)$ as well. Since $G(\alpha)$ is a compact set, there is $n_0\in\mathbb{N}$ such that $f^{n_0}\left(G(\alpha)\right)\subset f^{k_0}(U)$. It follows that
$f^{n_0+k}\left(G(\alpha)\right)\subset f^{k_0+k}(U)\subset f^{k_0}(U)\subset V_0(A)$ for every $k\in\mathbb{N}$. As a consequence, $\cup_{k\geq n_0}f^k\left(G(\alpha)\right)\subset V_0(A)$.
$\Box$

\section{Proof of Theorem \ref{thm:congugacy-semi-regular}}\label{s:proof-hom}

Suppose that homeomorphisms $f_1$, $f_2\in SsH(M^n)$ are conjugate. Since a conjugacy mapping $M^n\to M^n$ is a homeomorphism, the sets $A(f_1)$, $A(f_2)$, as well as the sets $R(f_1)$, $R(f_2)$ have the same dynamical embedding.

To prove the inverse assertion, let us suppose for definiteness that the basic sets $\alpha(f_1)$, $\alpha(f_2)$ are trivial while the sets $A(f_1)$, $A(f_2)$ have the same dynamical  embedding. Taking in mind that $A(f_1)$ and $A(f_2)$ are attracting sets, we see that there are neighborhoods $\delta_1$, $\delta_2$ of $A(f_1)$, $A(f_2)$ respectively, and a homeomorphism $h_0: \delta_1\to\delta_2$ such that
\begin{equation}\label{eq:conjugacy-in-neighborhoods-delta}
    h_0\circ f_1|_{\delta_1}=f_2\circ h_0|_{\delta_1}, \quad f_1(\delta_1)\subset\delta_1,\quad h_0(A(f_1))=A(f_2).
\end{equation}
Without loss of generality, one can assume that $\delta_1\subset B(A(f_1))$. Moreover, taking $\delta_1$ smaller if one needs, we can assume that $clos~\delta_1$ is a trapping region for $f_1$ of the set $A(f_1)$. By (\ref{eq:conjugacy-in-neighborhoods-delta}), one gets
 $$ f_2(clos~\delta_2)=f_2\circ h_0(clos~\delta_1)=h_0\circ f_1(clos~\delta_1)\subset h_0(\delta_1)=\delta_2. $$
Thus, $clos~\delta_2$ is a trapping region for $f_2$ of the set $A(f_2)$. As a consequence, we get the following generalization of (\ref{eq:conjugacy-in-neighborhoods-delta})
\begin{equation}\label{eq:iteration-conjugacy-in-neighborhoods-delta}
    h_0\circ f^k_1|_{\delta_1}=f^k_2\circ h_0|_{\delta_1},\quad k\in\mathbb{N}, \quad f_1(clos~\delta_1)\subset\delta_1,\quad h_0(A(f_1))=A(f_2).
\end{equation}

By Lemma \ref{lm:A-is-attractive-set-alpha-trivial}, there is the trapping region $T(\alpha_1)$ for $f^{-1}_1$ of the set $\widetilde{\alpha(f_1)}$ consisting of pairwise disjoint open $n$-balls $b_1$, $\ldots$, $b_r$ such that each $b_i$ contains a unique periodic point $q_i$ from $\alpha(f_1)$. In addition, the region $B(\widetilde{\alpha(f_1)})\setminus\widetilde{\alpha(f_1)}$ has the generating set $G(\alpha_1)$ consisting of pairwise disjoint closed $n$-annuluses $a_1$, $\ldots$, $a_r$ such that $a_i=clos~f^{p_i}_1(b_i)\setminus b_i$ where $p_i\in\mathbb{N}$ is a minimal period of the periodic point $q_i$.

Due to Lemma \ref{lm:neighbor-for-A-is-attractive-set-alpha-trivial}, one can assume without loss of generality that $G(\alpha_1)\stackrel{\rm def}{=}G_1\subset\delta_1$. Hence,
 $$ A(f_1)\bigcup\left(\cup_{k\geq 0}f^k(G_1)\right)=A(f_1)\bigcup N^+\subset\delta_1,\quad N^+=\cup_{k\geq 0}f^k(G_1). $$

According to Lemma \ref{lm:A-is-attractive-set-alpha-trivial}, $G_1$ is a generating set of the region $B(A(f_1))\setminus A(f_1)$. Let us show that $h_0(G_1)\stackrel{\rm def}{=}G_2$ is a generating set for the region $B(A(f_2))\setminus A(f_2)$. Take a point $z_2\in G_2$. There is a unique point $z_1\in G_1$ such that $h_0(z_1)=z_2$. Note that $z_2\notin A(f_2)$ since $z_1\notin A(f_1)$. Since $G_1\subset\left(B(A(f_1))\setminus A(f_1)\right)$, $f_1^k(z_1)\to A(f_1)$ as $k\to\infty$. It follows from (\ref{eq:iteration-conjugacy-in-neighborhoods-delta}) that
 $$ f_2^k(z_2)=f_2^k\circ h_0(z_1)=h_0\circ f_1^k(z_1)\to h_0(A(f_1))=A(f_2)\quad\mbox{ as }\quad k\to\infty . $$
Hence, $z_2\in B(A(f_2))$ and $G_2\in B(A(f_2))\setminus A(f_2)$.

Take an orbit $Orb_{f_2}\subset B(A(f_2))\setminus A(f_2)$. Since this orbit intersects a trapping region of $A(f_2)$, $Orb_{f_2}\cap\delta_2\neq\emptyset$. Therefore there exists a point $x_2\in Orb_{f_2}\cap\delta_2$. Since $h_0(A(f_1))=A(f_2)$ and $x_2\in B(A(f_2))\setminus A(f_2)$, the orbit $Orb_{f_1}$ of the point $x_1=h_0^{-1}(x_2)\subset\delta_1$ under $f_1$ belongs to $B(A(f_1))\setminus A(f_1)$. Hence, $Orb_{f_1}$ intersects $G_1$ at some point $w_1\in\delta_1$. Since $x_1$, $w_1\in Orb_{f_1}$, there is $k\in\mathbb{N}$ such that either $x_1=f_1^k(w_1)$ or $w_1=f_1^k(x_1)$. Suppose for definiteness that $w_1=f_1^k(x_1)$. Using (\ref{eq:conjugacy-in-neighborhoods-delta}), one gets
 $$ w_2=h_0(w_1)=h_0\circ f_1^k(x_1)=h_0\circ f_1^k\circ h_0^{-1}(x_2)=f_2^k(x_2)\in G_2\cap Orb_{f_2}. $$
Similarly one can prove that if $Orb_{f_2}$ intersects the interior of $G_2$, then $Orb_{f_2}$ intersects $G_2$ at a unique point, and if $Orb_{f_2}$ intersects the boundary of $G_2$ then $Orb_{f_2}$ intersects $G_2$ at two points. Thus, $G_2$ is a generating set for the region $B(A(f_2))\setminus A(f_2)$.

Set
 $$ \cup_{k\geq 0}f^{-k}_i(G_i)\stackrel{\rm def}{=}O^-(G_i),\quad \cup_{k\geq 0}f^{k}_i(G_i)\stackrel{\rm def}{=}O^+(G_i),\quad i=1,2. $$
We see that $O^-(G_i)\cup O^+(G_i)$ is invariant under $f_i$, $i=1,2$.
Given any point $x\in O^-(G_1)\cup O^+(G_1)$, there is $m\in\mathbb{Z}$ such that $x\in f^{-m}_1(G_1)$. Let us define the mapping
 $$ h: O^-(G_1)\cup O^+(G_1)\to O^-(G_2)\cup O^+(G_2) $$
as follows
 $$ h(x)=f^{-m}_2\circ h_0\circ f^m_1(x),\quad where\quad x\in f^{-m}_1(G_1). $$
Since $G_1$ and $G_2$ are generating sets, $h$ is correct. It is easy to check that
 $$ h\circ f_1|_{O^-(G_1)\cup O^+(G_1)}=f_2\circ h|_{O^-(G_1)\cup O^+(G_1)}. $$
It follows from (\ref{eq:conjugacy-in-neighborhoods-delta}) that
 $$ h: A(f_1)\cup O^-(G_1)\cup O^+(G_1)\to A(f_2)\cup O^-(G_2)\cup O^+(G_2) $$
is the homeomorphic extension of $h_0$ putting $h|_{A(f_1)}=h_0|_{A(f_1)}$. Moreover,
 $$ h\circ f^k_1|_{A(f_1)\cup O^-(G_1)\cup O^+(G_1)}=f^k_2\circ h|_{A(f_1)\cup O^-(G_1)\cup O^+(G_1)},\quad k\in\mathbb{Z}. $$

By Lemma \ref{lm:A-is-attractive-set-alpha-trivial}, $G_i$ is a generating set for the region
$B\left(\widetilde{\alpha(f_i)}\right)\setminus\widetilde{\alpha(f_i)} = B\left(A(f_i)\right)\setminus A(f_i)$ and $B\left(A(f_i)\right)\setminus A(f_i)=\cup_{k\in\mathbb{Z}}f^k_i(G_i)$, $i=1,2$. Thus, one gets the conjugacy $h: M^n\setminus\widetilde{\alpha(f_1)}\to M^n\setminus\widetilde{\alpha(f_2)}$ from $f_1|_{M^n\setminus\widetilde{\alpha(f_1)}}$ to  $f_2|_{M^n\setminus\widetilde{\alpha(f_2)}}$ :
\begin{equation}\label{eq:allmost-conjugacy}
    h\circ f^k_1|_{M^n\setminus\widetilde{\alpha(f_1)}}=f^k_2\circ h|_{M^n\setminus\widetilde{\alpha(f_1)}},\quad k\in\mathbb{Z}.
\end{equation}

Recall that the sets $\alpha(f_1)$, $\alpha(f_2)$ are periodic sources $\{\alpha_j(f_1)\}_{j=1}^{l_1}$, $\{\alpha_j(f_2)\}_{j=1}^{l_2}$ respectively. By Lemma \ref{lm:A-is-attractive-set-alpha-trivial}, the generating set $G_i$ consists of pairwise disjoint $n$-annuluses $a_j(f_i)$, $i=1,2$. Take an annulus $a_r(f_1)=a_r\subset G_1$ surrounding a source periodic point $\alpha_r(f_1))$ of minimal period $p_r$, $1\leq r\leq l_1$. Then the set $\bigcup_{k\geq 0}f_1^{-kp_r}(a_r)\cup\{\alpha_r(f_1))\}=D^n_r$ is a closed $n$-ball. Since
 $$ M^n\setminus B(A(f_2)) = M^n\setminus\left(A(f_2)\cup_{k\in\mathbb{Z}}f_2^k(G_2)\right) $$
consists of the source periodic points $\alpha(f_2)$, the annulus
 $$ \bigcup_{k\geq 0}f_2^{-kp_r}\circ h(a_r)=\bigcup_{k\geq 0}h\circ f_1^{-kp_r}(a_r) = D^*_r $$
surrounds a unique source periodic point $\alpha_{j(r)}(f_2)$ of the same minimal period $p_r$. Moreover, $D^*_r\cup\{\alpha_{j(r)}(f_2)\}$ is a closed n-ball. It implies the one-to-one correspondence $r\to j(r)$ inducing the one-to-one correspondence $j_0: \alpha_r(f_1))\to\alpha_{j(r)}(f_2))$. Since $\alpha_r(f_1))$ and $\alpha_{j(r)}(f_2))$ have the same period, one gets
\begin{equation}\label{eq:j-agree-with-f}
    j_0\left(f_1^k(\alpha_r(f_1)\right)=f_2^k\left(j_0(\alpha_r(f_1))\right)=f_2^k\left(\alpha_{j(r)}(f_2)\right),\quad 0\leq k\leq p_r.
\end{equation}

Put by definition, $h\left(\alpha_r(f_1)\right) = \alpha_{j(r)}(f_2))$. For sufficiently large $m\in\mathbb{N}$, the both $f^{-mp_r}_1(D^n_r)$ and $f^{-mp_r}_2(D^*_r)$ can be embedded in arbitrary small neighborhoods of $\alpha_r(f_1))$ and $\alpha_{j(r)}(f_2))$ respectively, because of $\widetilde{\alpha(f_1)}$ and $\widetilde{\alpha(f_2)}$ are repelling sets. Taking in mind (\ref{eq:j-agree-with-f}), it follows that $h: M^n\to M^n$ is a conjugacy from $f_1$ to $f_2$. This completes the proof.
$\Box$

\section{Some applications}\label{s:discus-and-appl}\nopagebreak

Following Smale \cite{Smale60a,Smale67}, we write $\sigma_1\succ\sigma_2$ provided $W^u(\sigma_1)\cap W^s(\sigma_2)\neq\emptyset$ where $\sigma_1$ and $\sigma_2$ are saddle periodic points.
Later on, we assume a surface $M^2$ to be closed and connected. Recall that a node is either a sink or a source.

\medskip
\textsl{Proof of Proposition \ref{prop:attractors-no-more-saddles-plus-one}} is by induction on $s_0$. First, we consider the case $s_0=0$. We have to prove that $k=1$. Suppose the contrary that is $k\geq 2$. According to \cite{GrinMedvPo-book,Plykin84} (see also \cite{GinesZhuzhoma1979,GinesZhuzhoma2005}), there are disjoint open sets $U_i$, $i=1,\ldots,k$, such that each $U_i$ is an attracting domain of $\Lambda_i$ with no trivial basic sets. Moreover, the boundary $\partial U_i$ consists of a finitely many simple closed curves. Therefore, $M^2\setminus\cup_{i=1}^kU_i$ is the disjoint union $\cup_{j\geq 1}K_j=G$ of compact connected sets $K_j$ where
$f^{-1}(G)\subset G$. Any iteration of $f$ has at least $k$ one-dimensional expanding attractors. Thus, without loss of generality, we can assume that
$f^{-1}(K_j)\subset K_j$ for every $K_j$. In addition, one can assume that any periodic isolated point is fixed and the restriction of $f$ on every invariant manifold of saddle isolated point preserves orientation.

Since $k\geq 2$ and $M^2$ is connected, there is a component of $G$, say $K_1$, and different sets $U_l$, $U_r$ such that $\partial K_1\cap\partial U_l\neq\emptyset$ and $\partial K_1\cap\partial U_r\neq\emptyset$ where
$U_l\cap U_r=\emptyset$. Any component of the boundary $\partial K_1$ is a circle. We see that there are at least two components of $\partial K_1$. Let us glue a disk to each boundary component of $\partial K_1$ to get a closed surface $\widetilde{K}_1$. Since $f^{-1}(K_1)\subset K_1$, one can extend $f|_{K_1}$ to an A-diffeomorphism $\widetilde{f}: \widetilde{K}_1\to\widetilde{K}_1$ with a unique sink in each disk we glued. Note that by construction, the non-wandering set $NW(\widetilde{f})$ of $\widetilde{f}$ consists of isolated nodal fixed points, and $NW(\widetilde{f})$ contains at least two sinks. According to \cite{Smale60a}, the surface $\widetilde{K}_1$ is the disjoint union of the stable manifolds of sinks and finitely many isolated sources (remark that the stable manifold of a source coincide with this source). This contradicts to the connectedness of $\widetilde{K}_1$ because of every stable manifold of a sink is homeomorphic to an open ball, and isolated sources do not separate the stable manifolds of two sinks. This contradiction proves that $k=1$ provided $s_0=0$.

Suppose the statement holds for $0,\ldots,s$ saddles. We have to prove this statement for $s_0=s+1$ saddles. Recall that due to \cite{Smale60a}, the isolated saddles endowed with the Smale partial order $\succ$. Since now the set of isolated saddles is not empty, there is a minimal saddle, say $\sigma$. Then the topological closure of $W^s(\sigma)$ is either a segment $I$ with the endpoints being two sources or a circle $S$ consisting of one source and $W^s(\sigma)$. In any cases, the both $I$ and $S$ are repelling sets. Let us consider this cases.

The segment $I$ has a neighborhood $U(I)=U$ homeomorphic to a disk such that $clos\,U\subset f(U)$. Note that $\sigma$ is inside of $U$. One can change $f$ inside of $U$ replacing $clos\,W^s(\sigma)$ by a unique source. One gets a diffeomorphism with $k$ expanding attractors and $s$ saddles. By the inductive assumption, $k\leq s+1\leq s_0<s_0+1$.

Similarly, the circle $S$ has a neighborhood $U(S)$ homeomorphic to an annulus such that $clos\,U(S)\subset f(U(S))$. Note that $\sigma$ belongs to $U(S)$. The manifold $M^2_1=M^2\setminus U(S)$ has two boundary components $M_1$, $M_2$ each homeomorphic to a circle. One can attach two disks $D^2_1$, $D^2_2$ along their boundaries to $M_1$, $M_2$ respectively to get a closed surface $\tilde{M}^2$. This surface either is connected or consists of two connected surfaces. Since $S$ is a repelling set, one can extend $f$ on $\tilde{M}^2$ to get a diffeomorphism $\tilde{f}: \tilde{M}^2\to \tilde{M}^2$ with $k$ expanding attractors and $s$ saddles. If $\tilde{M}^2$ is connected then the inductive assumption implies $k\leq s+1\leq s_0$. Let us consider the case when $\tilde{M}^2$ consists of two connected closed surfaces $\tilde{M}_1^2$, $\tilde{M}_2^2$. Suppose that $\tilde{M}_i^2$ contains $k_i$ expanding attractors and $s_i$ isolated saddles, $i=1,2$. Obviously, $k=k_1+k_2$ and $s_1+s_2=s$. By the inductive assumption, $k_i\leq s_i+1$, $i=1,2$. Hence, $k\leq (s_1+1)+(s_2+1)=s_1+s_2+2=s+2=s_0+1$. This concludes the proof.
$\Box$

\medskip
\textit{Proof of Theorem \ref{thm:many-expand-attr-one-less-saddle}.} Let us consider a structurally stable diffeomorphism $f: M^2\to M^2$ with the non-wandering set consisting of $k\geq 2$ one-dimensional expanding attractors $\Lambda_1$, $\ldots$, $\Lambda_k$, and isolated source periodic orbits, and $k-1$ saddle periodic points $\sigma_1$, $\ldots$, $\sigma_{k-1}$. Each $\Lambda_i$ has a neighborhood $U_i$ that is an attracting region of $\Lambda_i$. Then $M^2\setminus (\cup_{i=1}^{k-1}U_i)$ is the disjoint union $G=\cup_{j\geq 1}K_j$ of compact connected sets where $f^{-1}(G)\subset G$. Note that any positive iteration of $f$ has at least $k$ one-dimensional expanding attractors.  Obviously, any iteration of $f$ has the same number $k-1$ of saddle periodic points. Due to Proposition \ref{prop:attractors-no-more-saddles-plus-one}, any positive iteration of $f$ has no more than $k$  one-dimensional expanding attractors. Hence, any positive iteration of $f$ has exactly the same number $k$ of expanding attractors. This implies that every attractor $\Lambda_i$ is $C$-dense \cite{Anosov70,Robinson-book-99}. As a consequence, each unstable manifold $W^u(\cdot)\subset\Lambda_i$ is dense in $\Lambda_i$ \cite{Anosov70,GrinMedvPo-book}.

Take a connected component $K$ of the set $G$. The boundary $\partial K$ is the disjoint union of circles $c_1$, $\ldots$. By construction, this circles belong to the boundaries of the attracting regions $U_1$, $\ldots$, $U_k$. Therefore, one can glue a disk $d_j$ to each circle $c_j$ extending $f$ to $d_j$ with a sink inside of $d_j$. If $K$ is without an isolated saddles, then $K\cup_jd-J$ is a 2-sphere with a unique source and a unique sink \cite{GrinesGurevichZhPochinka2019}. Therefore if $K$ is without an isolated saddles, then $K$ is a disk with a unique source. Such a set $K$ we'll call a disk with no saddles. Now, take a neighborhood $U$ of some $\Lambda_i$. Suppose that all components of the boundary $\partial U$ attach to components of $G$ that are disks with no saddles. Then the union of $U$ and this disks gives a closed surface with  exactly one expanding attractor $\Lambda_i$. This contradicts to either the connectedness of $M^2$ or the inequality $k\geq 2$. Thus, given any neighborhood $U_i$ of $\Lambda_i$, the boundary $\partial U_i$ has a common part with the boundary $\partial K_j$ of some component $K\subset G$ which contains at least one isolated saddle.

Let $K$ be a component of $G$ containing a saddle $\sigma$ and $U$ a neighborhood of some $\Lambda_i$ such that $\partial K\cap\partial U\neq\emptyset$. Let us show that $W^u(\sigma)\cap W^s(\Lambda_i)\neq\emptyset$. Suppose the contrary. We know that $W^u(\sigma)\setminus\{\sigma\}$ belongs to stable manifolds of isolated periodic points lying in $K$. Then there is a saddle $\sigma_1\in K$ such that $\sigma_1\prec\sigma$, and the topological closure of $W^s(\sigma_1)$ is either a segment $I$ with the endpoints being two sources or a circle $S$ consisting of one source and $W^s(\sigma_1)$. In any cases, the both $I$ and $S$ are repelling sets. Therefore, $f$ can be changed inside of $K$ so that a diffeomorphism obtained has $k-2$ isolated saddles and $k$ one-dimensional expanding attractors. This contradicts Proposition \ref{prop:attractors-no-more-saddles-plus-one}. Thus,
$W^u(\sigma)\cap W^s(\Lambda_i)\neq\emptyset$.

Since $f$ is a structurally stable diffeomorphism, all intersections $W^u(\sigma)\cap W^s(x)$, $x\in\Lambda_i$, are transversal. It follows from $W^u(\sigma)\cap W^s(\Lambda_i)\neq\emptyset$ that there is $x\in\Lambda_i$ such that $W^u(\sigma)$ intersects transversally the stable manifold $W^s(x)$. Recall that the attractor $\Lambda_i$ is $C$-dense. Since any unstable manifold $W^u(\cdot)\subset\Lambda_i$ is dense in $\Lambda_i$, the topological closure of $W^u(\sigma)$ contains $\Lambda_i$, $clos\, W^u(\sigma)\supset\Lambda_i$.

Clearly that if $f_1$, $f_2$ are conjugate, then $\cup_{j=1}^{j=k-1}W^u(\sigma_j^{(1)})$, $\cup_{j=1}^{j=k-1}W^u(\sigma_j^{(2)})$ have the same dynamical embedding. Suppose that the sets $\cup_{j=1}^{j=k-1}W^u(\sigma_j^{(1)})$ and $\cup_{j=1}^{j=k-1}W^u(\sigma_j^{(2)})$ have the same dynamical embedding. It follows from above that
 $$clos\,\left(\cup_{j=1}^{j=k-1}W^u(\sigma_j^{(i)})\right)\supset\cup_{j=1}^{j=k}\Lambda_j^{(i)}, \quad i=1,2.$$
Since $A(f_i)=\cup_{j=1}^{j=k-1}W^u(\sigma_j^{(i)})\bigcup\left(\cup_{j=1}^{j=k}\Lambda_j^{(i)}\right)$, we see that
$$clos\, A(f_1)=clos\,\left(\cup_{j=1}^{j=k-1}W^u(\sigma_j^{(1)})\right),\quad clos\, A(f_2)=clos\,\left(\cup_{j=1}^{j=k-1}W^u(\sigma_j^{(2)})\right). $$
Therefore, the sets $A(f_1)$, $A(f_2)$ have the same dynamical embedding. As a consequence of Theorem \ref{thm:congugacy-semi-regular}, we have that $f_1$, $f_2$ are conjugate. This completes the proof.
$\Box$

\begin{figure}[h]
\centerline{\includegraphics[height=3.5cm]{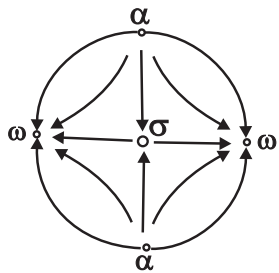}}
\caption{Phase portrait for $f\in MS(\mathbb{P}^2,1,1,1)$: the diametrically opposite points are identified.} \label{3-points}
\end{figure}

Consider $f\in MS(\mathbb{P}^2,1,1,1)$ with a unique saddle $\sigma(f)$. By definition, $f$ conjugates in some neighborhood of $\sigma(f)$ to a linear diffeomorphism with a saddle hyperbolic fixed point \cite{Pugh1962}. It easy to check that up to conjugacy there are exactly four such mappings :
 $$ T_1=
 \left\{\begin{array}{ccc}
         \bar{x} & = & \frac{1}{2}x \\
         \bar{y} & = & 2y,
       \end{array}\right.\qquad T_2=\left\{\begin{array}{ccc}
         \bar{x} & = & -\frac{1}{2}x \\
         \bar{y} & = & 2y,
       \end{array}\right.\qquad T_3=\left\{\begin{array}{ccc}
         \bar{x} & = & \frac{1}{2}x \\
         \bar{y} & = & -2y,
       \end{array}\right.\qquad T_4=\left\{\begin{array}{ccc}
         \bar{x} & = & -\frac{1}{2}x \\
         \bar{y} & = & -2y.
       \end{array}\right.
 $$
We'll say that the saddle $\sigma(f)$ is of the type $T_1$, $T_2$, $T_3$, $T_4$ respectively, see Fig.~\ref{3-points}.

\medskip
\textsl{Proof of Proposition \ref{prop:proj-plane-top-conj}}. Take $f\in MS(\mathbb{P}^2,1,1,1)$ with a unique saddle $\sigma(f)=\sigma$. The attracting set $A(f)$ is a closed curve consisting of an unstable manifold $W^u(\sigma)$ of a unique saddle $\sigma$ and a sink $\omega$. A neighborhood $U$ of $A(f)$ is homeomorphic to a M\"{o}bius band. Since $U$ contains only two fixed points, the saddle $\sigma$ and the sink $\omega$, the dynamics of $f|_{U}$ depends completely on a local dynamics of $f$ at $\sigma$ which is defined by one of the types $T_i$, $i=1,2,3,4$. Due to  Corollary \ref{cor:congugacy-regular}, diffeomorphisms $f_1$, $f_2\in MS(\mathbb{P}^2,1,1,1)$ are conjugate if and only if the types of their saddles coincide.

Choose any type $T_i\in\{T_1,T_2,T_3,T_4\}$. Let $B$ be a M\"{o}bius band with the middle closed curve $c_0$. There is a mapping $f_0: B\to B$ with the attracting set $c_0$ such that the non-wandering set of $f_0$ consists of a hyperbolic sink $\omega\in c_0$ and a hyperbolic saddle $\sigma\in c_0$ with $W^u(\sigma)=c_0\setminus\{\omega\}$. Note that the set $\mathbb{P}^2\setminus B$ is a 2-disk $D^2$. Since $c_0$ is an attracting set, one can extend $f_0$ to $f$ with a hyperbolic source in $D^2$. This gives $f\in MS(\mathbb{P}^2,1,1,1)$ desired.
$\Box$

\medskip
\textsl{Proof of Theorem \ref{thm:classif-M-2k-for-8-16}}.
1) Since $f$ has a unique saddle, the both $W^u_{\sigma_f}\cup\{\omega_f\}$ and $W^s_{\sigma_f}\cup\{\alpha_f\}$ are topologically embedded spheres denoted by $S^{k_1}$ and $S^{k_2}$ respectively. Due to \cite{MedvedevZhuzhoma2013-top-appl}, $k_1=k_2=k$, and the complements $M^{2k}(1,1,1)\setminus\left(W^u_{\sigma_f}\cup\{\omega_f\}\right)$, $M^{2k}(1,1,1)\setminus\left(W^s_{\sigma_f}\cup\{\alpha_f\}\right)$ homeomorphic to
an open $2k$-ball (see also, \cite{MedvedevZhuzhoma2016}). Thus, we have the embedding
 $$\varphi_u(f): S^k\to W^u_{\sigma_f}\cup\{\omega_f\}\subset M^{2k}(1,1,1),\quad \varphi_s(f): S^k\to W^s_{\sigma_f}\cup\{\alpha_f\}\subset M^{2k}(1,1,1). $$
Since the codimension of $S^k$ equals $k\geq 4$, $\varphi_u(f)(S^k)$ and $\varphi_s(f)(S^k)$ are locally flat spheres \cite{DavermanVenema-book-2009}. Hence, $\varphi_u(f)$ and $\varphi_s(f)$ are basic embedding.

2) According to Th\'{e}or\`{e}m d'approximation by Haefliger \cite{Haefliger1961}, we can assume without loss of generality that $\varphi(S^k)$ is a smoothly embedded $k$-sphere. Hence, there is a tubular neighborhood $T^{2k}$ of $\varphi(S^k)$ that is the total space of locally trivial fiber bundle $p: T^{2k}\to\varphi(S^k)$ with the base $S^k_0=\varphi(S^k)$ and a fiber $k$-disk $D^k$ \cite{Hirsch76}. Let $\vartheta_{ns}: S_0\to S_0$ be a Morse-Smale diffeomorphism with a unique sink $\omega_0$ and a unique source $N$, so-called a "north-south"\, diffeomorphism. The fiber $p^{-1}(N)$ is an open $k$-disk. Let $\psi_N: p^{-1}(N)\to p^{-1}(N)$ be the mapping with a unique hyperbolic sink at $N$ such that $clos\,\psi_N(p^{-1}(N))\subset\psi_N(p^{-1}(N))$ and $\cap_{j\geq 0}\psi^j_N(p^{-1}(N))=\{N\}$. Since $p$ is a locally trivial fiber bundle, one can extend $\psi_N$ and $\vartheta_{ns}$ to get the mapping $f_0: T^{2k}\to T^{2k}$ such that
a) $N$ is a hyperbolic saddle with $k$-dimensional local stable and unstable manifolds, and $\omega_0$ is a hyperbolic sink;
b) given any point $a\in T^{2k}\setminus p^{-1}(N)$, $f^l_0(a)$ tends to $\omega_0$ as $l\to\infty$; moreover, $S_0=\cap_{l\geq 0}f^l_0(T^{2k})$.

It was proved in \cite{MedvedevZhuzhoma2013-top-appl} that the boundary $\partial T^{2k}$ of $T^{2k}$ is a $(2k-1)$-sphere, say $S^{2k-1}$. Moreover, $S^{2k-1}$ bounds the ball $B^{2k}=M^{2k}(1,1,1)\setminus T^{2k}$. Take a point $a_0\in B^{2k}$. Since $B^{2k}$ is a ball, one can extend $f_0$ to $B^{2k}$ to get a mapping $f: M^{2k}(1,1,1)\to M^{2k}(1,1,1)$ with a unique hyperbolic source at $a_0$. It follows from (a) and (b) that we get the desired Morse-Smale diffeomorphism $f\in MS^{2k}(1,1,1)$ with the sink $\omega_0=\omega_f$, the saddle $N=\sigma_f$, and the source $a_0=\alpha_f$.

3) The last statement immediately follows from Corollary \ref{cor:congugacy-regular}.
$\Box$

\medskip
\textsl{Proof of Theorem \ref{thm:iterations-equiv-embed}}. Obviously, if $f_1^{k_1}$, $f_2^{k_2}$ are conjugate for some $k_1$, $k_2\in\mathbb{N}$, then the stable (unstable) manifolds $W^{s(u)}(\sigma_{f_1})$, $W^{s(u)}(\sigma_{f_2})$ have equivalent embedding. We have to prove the inverse assertion.

Suppose for definiteness that the unstable manifolds $W^u(\sigma_{f_1})$, $W^u(\sigma_{f_2})$ have equivalent embed\-ding. Let $T^{2k}_i$ be a tubular neighborhood of $S^k_i=W^s(\sigma_{f_1})\cup\{\omega_i\}$ that is the total space of locally trivial fiber bundle $p: T^{2k}_i\to S^k_i$ with the base $S^k_i$ and a fiber $k$-disk $D^k_i$, $i=1,2$. Here, $\omega_i$ is a unique sink of $f_i$, $i=1,2$. Note that $\sigma_{f_1}$, $\omega_i\subset S^k_i$ and $S^k_i$ is an attracting set of $f_i$, $i=1,2$. It follows that there is $k_i\in\mathbb{N}$ such that $f_i^{k_i}(T^{2k}_i)\subset T^{2k}_i$. Moreover, without loss of generality one can assume that the restrictions
 $$ f_i^{k_i}|_{W^s(\sigma_{f_i})}: W^s(\sigma_{f_i})\to W^s(\sigma_{f_i}),\quad f_i^{k_i}|_{W^s(\sigma_{f_i})}: W^u(\sigma_{f_i})\to W^u(\sigma_{f_i}),\quad i=1,2 $$
preserve orientation. Taking the neighborhood $T^{2k}_i$ smaller if necessary, one assume that $f_i^{k_i}$ near the saddle $\sigma_{f_i}$ conjugates a linear hyperbolic diffeomorphism due to the Grobman-Hartman Theorem \cite{Grobman1959,Grobman1962,Hartman1963} (see also \cite{Pugh1962}), $i=1,2$. Hence, $f_i^{k_i}$ is embedded into a flow near the saddle $\sigma_{f_i}$, $i=1,2$. Since $f_i^{k_i}(T^{2k}_i)\subset T^{2k}_i$, the both $f_1^{k_1}$ and $f_2^{k_2}$ are embedded into the flows, say $f^t_1$ and $f^t_2$, in the neighborhoods $T^{2k}_1$, $T^{2k}_2$ respectively. Clearly, the unstable manifolds $W^u(\sigma_{f_1})$, $W^u(\sigma_{f_2})$ are the unstable manifolds of $f^t_1$ and $f^t_2$. Since $W^u(\sigma_{f_1})$ and $W^u(\sigma_{f_2})$ have equivalent embedding, it follows from the proof of Theorem 2 \cite{MedvedevZhuzhoma2016} that the flows $f^t_1$ and $f^t_2$ are conjugate. This implies that $f_1^{k_1}$ and $f_2^{k_2}$ are conjugate.
$\Box$

%%%%%%%%%%%%%%%%%%%%%%%%%%%%%%%%%%%%%%%%%%%%%%%%%%%%%%%%%%%%%%%%%%%%%%%%%%%

\noindent
\textit{E-mail:} medvedev-1942@mail.ru

\noindent
\textit{E-mail:} zhuzhoma@mail.ru


\begin{thebibliography}{99}

\bibitem{Akin-book-1993}
\textbf{Akin E.} \textit{The General Topology of Dynamical Systems}. Graduate Studies in Math., v. \textbf{1}(1993), AMS.

\bibitem{AkinHurleyKennedy-book-2003}
\textbf{Akin E., Hurley M., Kennedy J.A.} \textit{Dynamics of Topologically Generic Homeomorphisms}.  Memoirs of AMS, v. \textbf{164}(2003), no. 783.

\bibitem{Anosov67}
\textbf{Anosov D.V.} Geodesic flows on closed Riemannian manifolds of negative curvature. {\it Trudy Mat. Inst. Steklov.} \textbf{90} (1967);
English transl., \textit{Proc. Steklov Inst. Math.} {\bf 90} (1969). %\MR{36:7157}
\bibitem{Anosov70}
\textbf{Anosov D.V.} On one class of invariant sets of smooth dynamical systems. \textit{Proc. Int. Conf. ``Nonlinear Oscillations", Qualitative Methods}, Kiev, Vol.~2, \textbf{1970}, 39-45.

\bibitem{AnosovZhuzhoma2005}
\textbf{Anosov D.V., Zhuzhoma E.} \textit{Nonlocal asymptotic behavior of curves and leaves of laminations on covering surfaces}. \textit{Proc. Steklov Inst. of Math.}, \textbf{249}(2005), p. 221.

\bibitem{ABZ}
\textbf{Aranson S., Belitsky G., Zhuzhoma E.} \textit{Introduction to Qualitative Theory of Dynamical Systems on Closed Surfaces}. Translations of Math. Monographs, Amer. Math. Soc., \textbf{153}(1996).

\bibitem{AransonPlykinZhirovZhuzhoma1997}
\textbf{Aranson S., Plykin R., Zhirov A., Zhuzhoma E.} Exact upper bounds for the number of one-dimensional basic sets of surface A-diffeomorphisms. \textit{Journ. Dynamical and Control Systems}, \textbf{3}(1997), no 1, 3-18.

\bibitem{JiagNiWang2004}
\textbf{Boju Jiang, Yi Ni, Shicheng Wang.} 3-manifolds that admit knotted solenoids as attractors. \textit{Trans. Amer. Math. Soc.}, \textbf{356}(2004), 4371-43-82.

\bibitem{Bowen1970b}
\textbf{Bowen R.} Topological entropy and Axiom A. \textit{Global Analysis, Proc. Sympos. Pure Math.}, Amer. Math. Soc., \textbf{14}(1970), 23-42.

\bibitem{DavermanVenema-book-2009}
\textbf{Daverman R.J., Venema G.A.} \textit{Embeddings in Manifolds}. GSM, Amer. Math. Soc., Providence, \textbf{106}(2009).

\bibitem{EellsKuiper62}
\textbf{Eells J., Kuiper N.} Manifolds which are like projective planes. \textit{Publ. Math. IHES}, \textbf{14}(1962), 5-46.

\bibitem{Fr70}
\textbf{Franks J.} Anosov diffeomorphisms. {\it Global Analisys. Proc. Symp. in Pure Math., AMS} {\bf 14}(1970), 61-94. %\MR{42:6871}

\bibitem{GrinMedvPo-book}
\textbf{Grines V., Medvedev T., Pochinka O.} \textit{Dynamical systems on 2- and 3-manifolds}. Dev. Math., v. 46, Springer, 2016, 295p.

\bibitem{GrinesGurevichZhPochinka2019}
\textbf{Grines V., Gurevich E., Pochinka O., Zhuzhoma E.} Classification of Morse-Smale systems and topological structure of the underlying manifolds. \textit{Russian Math. Surveys}, \textbf{74}(2019), no 1, 37-110.

\bibitem{GrinesMedvedevPochZh2015a}
\textbf{Grines V., Medvedev T., Pochinka O., Zhuzhoma E.} On heteroclinic separators of magnetic fields in electrically conducting fluids. \textit{Physica D: Nonlinear Phenomena}, \textbf{294}(2015), 1-5.

\bibitem{GinesZhuzhoma1979}
\textbf{Grines V., Zhuzhoma E..} The topological classification of orientable attractors on an $n$-dimensional torus. \textit{Russian Math. Surveys}, \textbf{34} (1978), no. 4, 163--164.
%\MR{80k:58071}
\bibitem{GinesZhuzhoma2005}
\textbf{Grines V., Zhuzhoma E.} On structurally stable diffeomorphisms with codimension one expanding attractors. \textit{Trans. Amer. Math. Soc.}, \textbf{357}(2005), 617-667.

\bibitem{Grobman1959}
\textbf{Grobman D.} On a homeomorphism of systems of differential equations. \textit{Dokl. Akad. Nauk, USSR}, \textbf{128}(1959), 5, 880-881.
\bibitem{Grobman1962}
\textbf{Grobman D.} Topological classification of neighborhoods of singular point in $n$-dimensional space. \textit{Matem. sbornik, USSR}, \textbf{56}(1962), 1, 77-94.

\bibitem{Haefliger1961}
\textbf{Haefliger A.} Plongements diff\'{e}rentiable de vari\'{e}t\'{e}s dans vari\'{e}t\'{e}s. \textit{Comment. Math. Helv.}, \textbf{36}(1961-1962), 47-82.

\bibitem{Hartman1963}
\textbf{Hartman P.} On the local linearization of differential equations. \textit{Proc. AMS}, \textbf{14}(1963), no 4, 568-573.

\bibitem{Hirsch76}
\textbf{Hirch M.} \textit{Differential Topology.} Sringer-Verlag, \textbf{1976}.

\bibitem{HirschPughShub77-book}
\textbf{Hirsch M., Pugh C., Shub M.} \textit{Invariant Manifolds}. Springer-Verlag (\textit{Lect. Nots Math.}), \textbf{1977}.

\bibitem{JimingMaBinYu2007}
\textbf{Jiming Ma, Bin Yu.} The realization of Smale solenoid type attractors in 3-manifolds. \textit{Topology and its Appl.}, \textbf{154}(2007), 3021-3031.

\bibitem{KuznetsovSP-2005}
\textbf{Kuznetsov S.P.} Examples of a physical system with a hyperbolic attractor of the Smale-Williams type. \textit{Phys. Rev. Lett.}, \textbf{95}(2005), 144101.

\bibitem{Leont-Maier37}
\textbf{Leontovich E.A., Maier A.G.} On trajectories defining a qualitative structure of decomposition of sphere into trajectories. {\it Doclady Acad. Sci. USSR}, {\bf 14}(1937), 5, 251-257.
\bibitem{Leont-Maier55}
\textbf{Leontovich E.A., Maier A.G.} On the scheme defining the topological structure of decomposition into trajectories. {\it Doclady Acad. Sci. USSR}, {\bf 103}(1955), 4, 557-560.

\bibitem{Mane88b}
\textbf{Ma\~n\'e R.} A proof of $C^1$ stability conjecture. \textit{Publ. Math. IHES} \textbf{66}(1988), 161-210.

\bibitem{Manning74}
\textbf{Manning A.} There are no new Anosov diffeomorphisms on tori. \textit{Amer. Journ. of Math.}, \textbf{96}(1974), 422-429.

\bibitem{MedvedevZhuzhoma2008-MIAN}
\textbf{Medvedev V., Zhuzhoma E.} Global dynamics of Morse-Smale systems. \textit{Proc. Steklov Inst. of Math.}, \textbf{261}(2008), 112-135.
\bibitem{MedvedevZhuzhoma2013-top-appl}
\textbf{Medvedev V., Zhuzhoma E.} Morse-Smale systems with few non-wandering points. \textit{Topology and its Applications}, \textbf{160}(2013), issue 3, 498-507.
\bibitem{MedvedevZhuzhoma2016}
\textbf{Medvedev V., Zhuzhoma E.} Continuous Morse-Smale flows with three equilibrium positions. \textit{Sbornik: Mathematics}, \textbf{207}(2016), 5, 702–723.

\bibitem{Milnor-1956}
\textbf{Milnor J.} On manifolds homeomorphic to the 7-sphere. \textit{Annals of Math.}, \textbf{64}(1956), no 2, 399-405.

\bibitem{Newhouse70}
\textbf{Newhouse S.} On codimension one Anosov diffeomorphisms. {\it Amer. J. of Math.}, \textbf{92}(1970), no 3, 761-770.

\bibitem{NikZ99}
\textbf{Nikolaev I., Zhuzhoma E.} \textit{Flows on 2-dimensional manifolds}. Lect. Notes in Math., {\bf 1705}, Springer \textbf{1999}.

\bibitem{Palis69}
\textbf{Palis J.} On Morse-Smale dynamical systems. \textit{Topology}, \textbf{8}(1969), no 4, 385-404.

\bibitem{Plykin84}
\textbf{Plykin R.V.} On the geometry of hyperbolic attractors of smooth cascades. \textit{Russian Math. Surveys}, \textbf{39}(1984), no 6, 85-131.

\bibitem{Po1886}
\textbf{Poincare H.} Sur les courbes d\'efinies par les equations differentielles. {\it J. Math. Pures Appl.}, {\bf 2}(1886), 151-217.

\bibitem{Pugh1962}
\textbf{Pugh C.} On a theorem of P.Hartman. \textit{Amer. Journ. Math.}, \textbf{91}(1962), no 2, 363-367.

\bibitem{PughWalkerWilson1977}
\textbf{Pugh C., Walker R.B., Wilson F.W.} On Morse-Smale approximations -- a counterexample. \textit{Journ. Diff. Equat.}, \textbf{23}(1977), 173-182.


\bibitem{Robinson-book-99}
\textbf{Robinson C.} \textit{Dynamical Systems: stability, symbolic dynamics, and chaos}. Studies in Adv. Math., Sec. edition, CRC Press, \textbf{1999}.

\bibitem{RobinsonWilliams71}
\textbf{Robinson C., Williams R.} Finite stability is not generic. {\it Dynamical Systems: Proc. Symp.} (Brazil, 1971), Academic Press, New York, London, 1973, 451-462.

\bibitem{Smale60a}
\textbf{Smale S.} Morse inequalities for a dynamical system. \textit{Bull. Amer. Math. Soc.}, \textbf{66}(1960), 43-49.
\bibitem{Smale67}
\textbf{Smale S.} Differentiable dynamical systems. \textit{Bull. Amer. Math. Soc.}, \textbf{73}(1967), 747-817.

\bibitem{Strogatz-book-1994}
\textbf{Strogatz S.} \textit{Nonlinear Dynamics and Chaos with Applications to Physixs, Biology, Chemistry and Engineering}. Persius Books, Cambridge, Massachusetts.

\bibitem{Williams1970a}
\textbf{Williams R.} The DA maps of Smale and structural stability. \textit{Global Anal., Proc. Symp. Pure Math.}, AMS, \textbf{14}(1970), 329-334.
\bibitem{Williams1974}
\textbf{Williams R.} Expanding attractors. \textit{Publ. Math. I.H.E.S.}, \textbf{43}(1974), 169-203.

\bibitem{ZhIsaenkova2011-transl}
\textbf{Zhuzhoma E.V., Isaenkova N.V.} Zero-dimensional solenoidal base sets. \textit{Sbornik: Mathematics}, \textbf{202}(2011), no 3, 351-372.

\end{thebibliography}
\end{document}